\documentclass{article}
\usepackage{epsfig}
\usepackage{amscd}
\usepackage{latexsym}
\usepackage{amsmath}
\usepackage{amssymb}
\usepackage{theorem}
\usepackage{varioref}

\usepackage{subfigure}

\newtheorem{theorem}{Theorem}[section] 
\newtheorem{lemma}[theorem]{Lemma}

\newtheorem{corollary}[theorem]{Corollary}
\newtheorem{conjecture}[theorem]{Conjecture}

{\theorembodyfont{\rmfamily}
\theoremstyle{plain}
\newtheorem{definition}[theorem]{Definition}
\newtheorem{example}[theorem]{Example} 
\newtheorem{problem}[theorem]{Problem} 
\newtheorem{remark}[theorem]{Remark}
\newtheorem{ack}{Acknowledgements}
}

\newcommand{\G}{\ensuremath{\mathcal G}}
\newcommand{\A}{\ensuremath{\mathcal A}}
\newcommand{\B}{\ensuremath{\mathcal B}}
\renewcommand{\S}{\ensuremath{\mathcal S}}
\renewcommand{\to}{\longrightarrow}
\newcommand{\C}{\ensuremath{\mathbb C\,}}
\newcommand{\Z}{\ensuremath{\mathbb Z}}
\renewcommand{\P}{\ensuremath{\mathbb P}}
\newcommand{\K}{\ensuremath{\mathbb K}}
\newcommand{\Q}{\ensuremath{\mathbb Q}}
\newcommand{\codim}{\operatorname{codim}}
\newcommand{\edge}{\operatorname{edge}}
\newcommand{\lop}{\operatorname{loop}}

\newcommand{\we}{\ensuremath{\wedge}} 
\renewcommand{\bar}[1]{\ensuremath{\overline{#1}}}

\newcommand{\qed}{\hfill \mbox{$\Box$}\medskip}
\newenvironment{proof}{\noindent {\it proof:}}{\qed \par}

\title{Parallel connections and bundles of 
arrangements \thanks{revised July 17, 2000.}}

\author{Michael J. Falk \and Nicholas J. Proudfoot\thanks{research conducted in part under an NSF {\em
Research Experiences for Undergraduates} grant.}}

\date{}

\begin{document}
 
\maketitle

\begin{abstract}
Let \A\ be a complex hyperplane arrangement, and let $X$ be a modular 
element of arbitrary rank in the intersection lattice of \A.
Projection along $X$ restricts to a fiber bundle 
projection 
of the complement of \A\ to the 
complement of the localization $\A_X$ of \A\ at $X$. We identify the fiber 
as the 
decone of a realization of the complete principal truncation of the 
underlying matroid of \A\ along the flat corresponding to $X$. We 
then generalize to this setting several properties of strictly linear fibrations, the case 
in which $X$ has corank one, including the triviality of the monodromy 
action on the cohomology of the fiber.
This gives a topological 
realization of results of Stanley, Brylawsky, and Terao on modular 
factorization. We also show that 
(generalized) parallel connection of matroids corresponds to pullback 
of 
fiber bundles, clarifying the notion that all examples of 
diffeomorphisms of complements of inequivalent arrangements result from
 the triviality of the restriction of the Hopf 
bundle to the complement of a hyperplane. 
The modular fibration theorem also yields a new method for identifying 
$K(\pi,1)$ arrangements 
of rank greater than three. We exhibit new families of $K(\pi,1)$ arrangements, 
providing more evidence for the 
conjecture that factored arrangements of arbitrary rank are 
$K(\pi,1)$.
\end{abstract}

\begin{section}{Introduction}

Let $V$ be a vector space over a field \K. An {\em arrangement} \A\ 
in $V$ is a finite collection of linear hyperplanes in $V$. The {\em 
complement} $M=M(\A)$ of \A\ is $V - \bigcup \A$. A set of 
hyperplanes \B\ is {\em dependent} if the $\codim(\bigcap 
\B)<|\B|$. These dependent sets determine a matroid $G(\A)$ with ground set \A, 
the {\em underlying matroid} of \A. Alternatively, $G(\A)$ is
the linear matroid realized by the projective point configuration 
$\A^*$ in $\P(V^*)$ determined by the defining linear forms for the hyperplanes	 
of \A. 

In case $\K=\C$ the complement $M(\A)$ is a connected manifold whose 
topology has been studied in great detail. In this case 
there is a strong connection between the topological structure of 
$M(\A)$ 
and the underlying matroid $G(\A)$. The paradigmatic result along 
these 
lines is that the cohomology of $M(\A)$ has a presentation depending 
only on $G(\A)$, 
with the consequence that the Poincar\'e series of the cohomology 
ring of 
$M(\A)$ essentially coincides with the characteristic polynomial of 
$G(\A)$ 
\cite{OS1}. It has become clear that 
techniques and constructions from matroid theory can have interesting 
and surprising implications for the topology of hyperplane 
complements. In this paper we interpret 
the matroidal notions of modular flat, principal truncation, and 
generalized parallel connection in 
this vein, in terms of 
bundles 
of complex hyperplane arrangements, their fibers, and pullbacks via 
inclusion maps.

\medskip
Henceforth we restrict our study to complex arrangements.
The {\em intersection lattice} $L=L(\A)$ of \A\ is the set of 
subspaces $X$ of $\C^\ell$ which are intersections of hyperplanes of \A, 
$X=\bigcap \B$ for $\B\subseteq \A$, partially 
ordered by reverse inclusion. The 
smallest element of $L$ is $O_{L}=\C^{\ell}$, the empty intersection, 
and the largest element of $L$ is $1_{L}=\bigcap \A$.
For $X,Y \in L(\A),$ the {\em join} $X\vee 
Y$ is $X\cap Y$ and the {\em meet} $X\we Y$ is $\bigcap\{H\in \A \ | 
\ 
H\supseteq X+Y\}$. The 
rank function $r$ of $L$ is given by 
$r(X)=\codim(X)$, and
the semimodular law holds: $$r(X\we 
Y)+r(X\vee Y)\leq r(X)+r(Y)$$ for $X,Y\in L$. 
Then $L$ is a geometric lattice, isomorphic to the lattice of 
flats of the matroid $G(\A)$, via the identification of		
$X\in L$ with the flat $$\A_{X}=\{H\in \A \ | \ H\supseteq X\}.$$	We will	
often refer to elements $X\in L(\A)$ as flats, tacitly identifying $X$ 
with $\A_X$. For instance, ``point'' and ``line'' refer to flats of rank 
one and two. The {\em corank} of a flat $X$ is $r(1_{L})- r(X)$, and ``copoints'' 
and ``colines''  are flats of corank one and two.

When equality holds in the formula above, $(X,Y)$ is called a {\em 
modular 
pair}. An element $X\in L(\A)$ is {\em modular} if $(X,Y)$ is a 
modular pair for every $Y\in L(\A)$. This is equivalent to the 
condition that $X+Y$ be an element of $L$ for every $Y\in L$. 
Let $\pi$ be the linear projection of $\C^{\ell}$ 
onto the quotient $\C^{\ell}/X$. 
Modularity of 
$X$ implies that fibers of $\pi$, the parallel translates of $X$, 
intersect each $Y\in L(\A)$ in the same way, independent of position. 
This observation was already made by Terao in \cite{T2}, who proved 
that $\pi|_{M(\A)}$ is a fiber bundle projection in case $X$ has 
corank one. 
But, in fact, it is easy to show that modularity 
of $X$ is equivalent to $\pi$ being a map of stratified spaces, 
under the natural stratifications of $\C^{\ell}$ and $\C^\ell/X$ 
determined by $\A$ and $\A_{X}$. Being a linear projection, it is 
trivial to show $\pi$ restricts to a submersion on each stratum. 
L.~Paris showed how to extend $\pi$ to a proper map of stratified 
spaces. Then  
Thom's Isotopy Lemma implies that $\pi|_{M(\A)}$ is a fiber 
bundle projection for $X$ a modular flat of arbitrary rank. 

This fibration result interpolates between two well-known extreme 
cases. In case $X$ is a modular copoint, the result was proven in 
\cite{T2}, as already mentioned. This case gives rise to the notion 
of 
supersolvable arrangement, and its connection with fiber-type 
arrangements \cite{FR2}, a much-studied class \cite{Stan,JT,BZ2,CS5}. 
In case $X$ is a point, i.e., a hyperplane of \A, then $X$ is 
automatically modular, and the fibration is just the restriction of 
the defining form $\phi: \C^\ell \to \C$ of the hyperplane $X$. This 
gives rise to the well-known elementary ``cone-decone'' construction 
\cite{OT}. The restriction $\phi|_{M(\A)}: M(\A) \to \C^*$ is in fact a trivial 
fibration, with fiber  
isomorphic to the complement in $\C^{\ell-1}$ of an affine 
arrangement, the {\em decone} of \A. 

The general modular fibration theorem was proved by L.~Paris \cite{P3}.
At the same time, we were independently conducting the research 
reported on in this paper \cite{Proud}, and had arrived at the same 
conclusion, only 
to later discover an error in our treatment of the proper extension 
of $\pi$. We sketch the argument here, and refer the reader to \cite{P3} 
for a complete proof, concentrating instead on other 
structural results and consequences of the theorem.

In \cite{T2}, Terao establishes the result for 
modular copoints, and proves that for general modular $X$ the fibers 
of $\pi$ have the same combinatorial type. But
he specifically remarks that a proof of local triviality in the general 
case is not at hand. See Remark \ref{terao}.
The proof of Corollary~\ref{bundle} is in a sense a parametrized 
version of the argument of  \cite{R1}, where 
stratification techniques were first used in the theory of 
arrangements, several years after Terao's work.  

The characteristic polynomial of a lattice was defined by G.-C.~Rota. 
The characteristic polynomial of a matroid is the 
characteristic polynomial of its lattice of flats.
The modular flat $X$ gives rise to a factorization of the 
characteristic polynomial of $G(\A)$ over the integers, with one 
factor given by the characteristic polynomial of $G(\A_X)$. This is 
Stanley's modular factorization theorem \cite{Stan3}. Brylawski \cite{Bry2}
identified the other factor as the characteristic polynomial of a 
related matroid, the {\em complete principal truncation} $\bar{T}_X(G)$ 
\cite[Section 7.4]{Wh1} of 
$G=G(\A)$ along $X$, divided by $(t-1)$. The complete principal 
truncation 
is obtained by successively adjoining generic points on the specified 
flat 
and contracting on the new points. Technically this is a matroid with multiple 
points; when we refer to $\bar{T}_X(G)$ we will always mean the 
associated simple matroid (with the same lattice of flats).

We show in Theorem~\ref{trunc}
that the fiber of the bundle map $\pi|_{M(\A)}$ is the 
complement of the decone of an arrangement realizing the complete 
principal truncation of $G(\A)$ on the flat $X$. In addition, just 
as in the corank-one case \cite{FR2}, the monodromy of the bundle is 
shown to act trivially 
on the cohomology of the fiber (Theorem~\ref{trivial}). Then the $E_{2}$ term in the Leray-Serre 
spectral sequence of $\pi|_{M(\A)}$ is isomorphic to the tensor 
product of the cohomology of the base $M(\A_{X})$ with that of the fiber. 
Using the identity relating 
characteristic polynomials and Poincar\'e polynomials, we obtain 
a topological 
interpretation of the Stanley-Brylawski and Terao factorization results.
In fact, the factorization of the characteristic polynomial implies 
that the spectral sequence degenerates at the $E_{2}$ term, just as 
in the corank-one case, although we have no topological proof of this 
fact (Remark \ref{gysin}).

In the corank-one situation, the monodromy of the bundle gives rise to
a ``braid monodromy'' homomorphism from $\pi_1(M)$ to the (pure) braid 
group on $n$ strands, where $n=|\A-\A_X|$. In the 
general case the analogue of this braid  
monodromy takes values in the fundamental group of the matroid 
stratum of the Grassmannian, or equivalently, the projective 
realization space, 
of the complete principal truncation $\bar{T}_X(G)$. See Remark 
\ref{monodromy}.
 
The current research grew out of an attempt to clarify and generalize 
the construction of \cite{EF}, which involved arrangements whose 
matroids are parallel connections. We began by studying the matroidal 
notion of {\em generalized parallel connection}. Loosely speaking, 
this is the free sum of two matroids along a common flat. This free 
sum 
is well-defined if and only if the flat is modular in one of the matroids. 
Thus we were led to the consideration of modular flats. The 
combinatorial study of Sections \ref{proj} and \ref{parallel}
formed the main part of an NSF {\em 
Research Experiences for Undergraduates} project in the summer of 1997. 
This work was reported on in \cite{Proud}, which provided the groundwork 
for this paper. 

Given the modular fibration result, we show that generalized 
parallel connection, in a natural 
realization in terms of complex arrangements, corresponds to the 
pullback of fiber bundles (Theorem \ref{pullback}). The construction 
of \cite{EF}, which yields diffeomorphisms of the complements of 
arrangements with non-isomorphic matroids, uses ordinary parallel 
connection, in which the identified flats are points. Then the 
diffeomorphisms of \cite{EF} are a consequence of two elementary 
observations, that the cone-decone construction yields a trivial 
bundle, and that the pullback of a trivial bundle is trivial.

\medskip
When the base and fiber of a modular fibration are both 
aspherical, it follows that the complement $M(\A)$ is also 
aspherical. 
In this case \A\ is called a $K(\pi,1)$ arrangement. The problem of 
identifying $K(\pi,1)$ arrangements has been an important one in the 
study of complex arrangements. There are two well-known classes of 
$K(\pi,1)$ arrangements, the supersolvable ones, which are abundant 
in all ranks, and the simplicial ones, which are rare in ranks 
greater than three. Other techniques for identifying 
$K(\pi,1)$ arrangements are mostly restricted to arrangements of rank three.
See \cite{F1,FR1,FR4} for further exposition of the $K(\pi,1)$ problem.

Corollary~\ref{bundle} provides a method for identifying $K(\pi,1)$ 
arrangements in ranks greater than three. In the final section we 
exhibit new families of such arrangements, arising from the work of 
P.~Edelman and V.~Reiner \cite{ER} and D.~Bailey \cite{Bail} on 
threshold graphs and subarrangements of the Coxeter arrangement of 
type $B_\ell$. By the classification result of Bailey, these 
new examples are all ``factored'' \cite[Section 3.3]{OT}. So our result 
provides more evidence for the conjecture that factored arrangements of arbitrary rank 
are $K(\pi,1)$ \cite{FR4}.  We also give an example of an arrangement 
of 
rank four which has two different modular colines. Then 
Corollary~\ref{bundle} 
implies that a certain arrangement of rank three, the cone of one 
of the fibers, is not $K(\pi,1)$, an arrangement to which existing 
techniques do not apply.

The search for examples of high-rank $K(\pi,1)$ arrangements was 
motivated by a suggestion of G.~Ziegler several years ago concerning 
counterexamples to the ``homotopy type conjecture,'' that complex arrangements 
with the same underlying matroid should have homotopy-equivalent 
complements \cite{OT}. This idea, laid out in \cite{Proud}, is to find 
$K(\pi,1)$ arrangements of high rank, whose underlying matroids have 
different characteristic polynomials, but have isomorphic generic 
rank-three truncations. Then generic 3-dimensional sections of these 
arrangements will have isomorphic underlying matroids, but 
non-isomorphic fundamental groups. Unfortunately we are so far unable to 
construct such examples using the technique of this paper.

\end{section}

\begin{section}{Projections and principal truncations}
\label{proj}
In this section we establish terminology and 
analyze the combinatorics associated with projections of hyperplane 
arrangements. Let $\A$ be a 
central arrangement of hyperplanes in $\C^\ell$. Let $L=L(\A)$ be the 
intersection lattice of \A, consisting of subspaces of $\C^\ell$ as 
described in the introduction. 

Let $X\in L$. Let $\A_{X}=\{H\in \A \ | \ H\supseteq X\}$, and 
let $\pi: \C^\ell \to \C^\ell/X$ be the natural 
projection. Note that $\pi$  maps each hyperplane $H\in \A_{X}$ to a 
hyperplane 
of $\C^{\ell}/X$. Henceforth we consider the arrangement 
$\A_{X}$ to be an arrangement in $\C^{\ell}/X$.

We shall have occasion to study arrangements formed by the 
intersections of the hyperplanes of $\A$ with a given affine subspace $S$. This 
induced arrangement in $S$ is called the {\em restriction} of \A\ to 
$S$, denoted $\A^S$. 

We start by describing the combinatorial structure of the 
affine arrangement $\A_{\pi}$ 
formed by restricting \A\ to a generic fiber of $\pi$. 
This requires 
some discussion of the cone-decone construction of \cite{OT}, and a 
description of the matroid construction of principal truncation along 
a flat.

There is a natural correspondence between arrangements of linear 
hyperplanes in $\C^{\ell}$ and arrangements of affine hyperplanes in 
$\C^{\ell-1}$. The analytic operations need not concern us here; they 
are described in detail in \cite{OT} and \cite{EF}. 
One places a copy of a given affine 
$(\ell-1)$-arrangement \A\ into $\{1\} \times \C^{\ell-1}\subseteq 
\C^{\ell}$. Then replace each of the affine subspaces in this copy of 
\B\  by its linear span in $\C^{\ell}$
(i.e., ``cone over the origin'') and adjoin the 
``hyperplane at infinity'' $\{0\}\times \C^{\ell-1}$, to 
obtain a central arrangement $c\A$, the {\em cone} of \A, in 
$\C^{\ell}$. 
The inverse operation, called ``deconing,'' takes a central 
arrangement 
\A\ to its projective image, and dehomogenizes relative to a 
hyperplane $H_{\infty}\in \A$ to obtain an affine 
$(\ell-1)$-arrangement $d\A$.

The intersections of hyperplanes of $d\A$ form a geometric 
semilattice \cite{WW}, isomorphic to a subposet of $L(\A)$.
Specifically, $$L(d\A)\cong 
\{Y\in L(\A) \ | \ Y\not\geq H_{\infty}\}.$$

The fiber arrangement $\A_{\pi}$ is an affine arrangement of 
dimension $\ell-r(X)$. We will show that the underlying matroid of 
the 
cone $c\A_{\pi}$ is the complete principal truncation of the matroid 
$G(\A)$ 
along the flat $\A_{X}$. 

The {\em principal truncation} $T_{F}(G)$ of a matroid $G$ along a 
flat $F$ is 
constructed by adding a generic point $p$ on the flat $F$ and then 
contracting $G$ on $p$ \cite[Section 7.4]{Wh1}. The result may be a 
matroid with multiple points. We tacitly simplify the resulting 
matroid, by removing any multiple points. This does not affect the 
intersection lattice, characteristic polynomial, or Orlik-Solomon 
algebra.

This operation can be 
iterated. The {\em 
complete} principal truncation ${\bar T}_{F}(G)$ is the result of 
$r(F)-1$ successive principal truncations on $F$, so that $F$ reduces 
to a point. Equivalently, one can add $r(F)-1$ generic points to $F$ 
and contract $G$ on the flat spanned by the new points. Contraction 
of a matroid on a point corresponds to projection of a projective 
point configuration from one of its points, or restriction of a 
hyperplane arrangement to one of its hyperplanes. 

\begin{theorem} Let $X\in L$ and let $\A_{\pi}$ be the affine 
arrangement obtained by restricting \A\ to a generic fiber of 
$\pi: \C^{\ell} \to \C^\ell/X$. Then the matroid $G(c\A_{\pi})$ is 
isomorphic to the complete principal truncation ${\bar T}_{X}(G)$ of 
$G=G(\A)$ along the flat $\A_{X}$.
\label{trunc}
\end{theorem}

\begin{proof}
Dualizing the description of complete principal truncation to 
hyperplane 
arrangements, we see that ${\bar T}_{X}(G)$ is the matroid of the 
arrangement $\A^P$ obtained by choosing a generic subspace $P$ of 
codimension $r(X)-1$ containing $X$, and restricting \A\ to $P$. Then 
$P$ has dimension $\dim(X)+1$, $X\cap P$ is a hyperplane of 
$\A^P$, and an affine translate of $X\cap P$ is a generic fiber 
of $\pi$. It follows that $d(\A^P) \cong \A_\pi$, so $\A^P 
\cong 
c\A_{\pi}$.
\end{proof}

\begin{definition} A pair $(X,Y)$ forms a {\em modular pair} in $L$ 
if $$r(X\vee Y)+r(X\we Y) = r(X)+r(Y).$$ An element $X\in L$ is {\em 
modular} if $(X,Y)$ is a modular pair for every $Y\in L$.
\end{definition}

The following lemma is the key to the proof of the modular fibration theorem, 
and is trivial to prove.

\begin{lemma} Let $X,Y \in L$. Then $(X,Y)$ is a modular pair if and 
only if $X+Y\in L$. \qed
\label{sum}
\end{lemma}

When $X$ is modular, the  conclusion of Theorem \ref{trunc} holds for 
every fiber of $\pi$ over points not in $\bigcup \A_{X}$.
To prove this we need to describe the rank function $r_{T}$ on the 
lattice of 
flats $L({\bar T}_{X}(G))$.
According to \cite[Proposition 7.4.9]{Wh1}, 
the set $L({\bar T}_X(G))$ can be identified with 
$\{Y\in L \ | \ X\we Y = 0_{L(\A)}
\ \text{or} \  Y\geq X\}$. With this identification the 
rank function $r_{T}$ is given by $$r_T(Y) := \begin{cases} r(Y) \quad 
\text{ if} \
X\we Y = 0_{L(\A)}, \\  r(Y) - r(X) + 1 \quad \text{if} \  Y\geq X. 
\end{cases}$$

\begin{theorem} Suppose $X$ be a modular flat. Let $\bar{v}=v+X \in 
(\C^{\ell}/X) 
-\bigcup \A_{X}$ and let $\A_{\bar{v}}$ be the restriction of \A\ to 
$\pi^{-1}(\bar{v})$. Then the intersection lattice $L(c\A_{\bar{v}})$ 
is isomorphic to 
$L({\bar T}_X(M)).$
\label{lattice} 
\end{theorem} 

\begin{proof} As in the proof of Theorem \ref{trunc}, the
arrangement $c\A_{\bar{v}}$ can be identified with the restriction 
$\A^P$ of \A\ to the linear subspace $P$ of codimension $r(X)-1$ spanned 
by $X$ and $v$. Then $L(c\A_{\bar{v}})=\{P\cap Y \ | \ Y\in L(\A)\}.$ 
There are three cases. 
\begin{list}{}
\item{Case 1.} Suppose $Y\in L(\A)$ satisfies $X\we Y= 
0_{L(\A)}$. By modularity of $X,$ 
$X+Y=\C^{\ell}$. Then there exists $y\in Y$ such that 
$\bar{v}=\bar{y}=y+X$. Then $P\cap Y=(\C v+X)\cap Y=(\C y+X)\cap 
Y=\C y+(X\cap Y)$. Since $\bar{v}\not \in \bigcup\A_{X}$, $y\not \in 
X$, so 
\begin{equation*}
\begin{split}
\codim_{P}(P\cap Y)&=\codim_{\C y+X}(\C y+(X\cap Y)) \\
&=\codim_{X}(X\cap Y)\\
&=\codim_{\C^{\ell}}(X \cap  Y)-\codim_{\C^{\ell}}(X)\\
&=r(X\vee Y)-r(X)\\
&=r(Y),\\
\end{split}
\end{equation*}
the last equality by modularity of $X$.
\item{Case 2.} Suppose $Y\geq X$. Then $Y\subseteq X\subseteq P$ so
\begin{equation*}
\begin{split}
\codim_{P}(P\cap Y)&=\codim_{P}(Y)\\
&=\codim_{X}(Y)+\codim_{P}(X)\\
&=r(Y)-r(X)+1.\\
\end{split}
\end{equation*}
\item{Case 3.} Suppose $0_{L(\A)}<X\we Y<X$. Then $X\cup Y\subseteq 
H$ for 
some $H\in \A$. Note that $v\not \in H$, since $\bar{v}\not\in 
\bigcup \A_{X}$. It follows that $P\cap H=(\C v+X)\cap H=X$ since $X\subseteq 
H$ while $v\not\in H$. Since $P\cap Y\subseteq P\cap H$ we have 
$P\cap Y=P\cap Y'$ for $Y'=X\vee Y\geq X$, which case is treated above.
\end{list}
These calculations verify that $L(c\A_{v})$ can be identified with 
$L(\bar{T}_X(G))$ as described above, with the same rank function.
\end{proof}

\begin{remark} The same calculations as in case 1 above can be used 
to show that 
the converse of Theorem \ref{lattice} also holds. That is, $X$ is 
modular if the lattice $L(c\A_{\bar{v}})$ is constant over $M(\A_X)$.
This will be used to identify modular flats in the examples of 
Section~\ref{examples}.\qed
\label{converse}
\end{remark}

Let $M(\A)=\C^{\ell}-\bigcup \A$ and $M(\A_{X})=(\C^\ell /X) - \bigcup 
\A_{X}$. 
Note that $\pi$ maps $M(\A)$ onto $M(\A_{X})$.
\begin{corollary} The fibers of $\pi|_{M(\A)}: M(\A) \to M(\A_{X})$ 
are 
diffeomorphic.
\label{diffeo}
\end{corollary}

\begin{proof} The fiber of $\pi|_{M(\A)}$ over $\bar{v}$ is the 
complement of the arrangement $\A_{\bar{v}}$ in $\pi^{-1}(\bar{v})\cong 
\C^{\ell-r(X)}$. Since the base $M(\A_{X})$ is path-connected, 
Theorem \ref{lattice} implies 
that the arrangements $c\A_{v}$ are lattice-isotopic. Then the 
assertion 
follow from \cite{R1}.
\end{proof}

\begin{remark} Theorem \ref{lattice} was essentially proved by Terao 
in 
\cite{T2}. Our result explicitly identifies 
the lattice. In case $X$ is a copoint, Corollary \ref{diffeo} follows 
without 
using Randell's lattice isotopy theorem, which had not been 
discovered at 
the time of Terao's work. In fact Corollary \ref{diffeo} and the 
fibration result Corollary \ref{bundle} of the next section confirm the suggestion stated 
after 
Proposition 2.12 of \cite{T2}. The proof of Corollary \ref{bundle} 
uses the stratification technique first introduced to 
arrangement theory by Randell in his proof of the isotopy theorem.
\label{terao}
\end{remark}

\end{section}

\begin{section}{Modular flats and fibrations}
\label{fibration}
The arrangement \A\ defines a stratification \S\ of $\C^{\ell}$:
$$\C^{\ell}=\bigcup \{S_Y \ | \ Y\in L\}$$
 whose strata $S_Y$ are given by $$S_Y=Y -\bigcup_{Z>Y} 
Z.$$
Thus $S_Y$ is a connected dense open subset of the linear space $Y$. 
In particular, $S_{Y}$ is a smooth submanifold of $\C^{\ell}$. 
Note that the 
closed stratum $\bar{S_{Y}}$ is equal to $Y$.
Also $S_{Y}\cap \bar{S_{Z}}\not = \emptyset$ if and only if 
$S_{Y}\subseteq \bar{S_{Z}}$ if and only if $Y\geq Z$.

This stratification satisfies Whitney's conditions (a) and (b) 
\cite{GorMac}. Indeed these conditions involve tangent and 
secant lines, and tangent spaces to strata, which are trivial to 
verify because $S_Y$, as an open subset of the linear space $Y$, has 
tangent space at any point equal to
$Y$.

Let $X\in L$ be a modular element of rank $p$. 
We may identify $\C^{\ell}/X$ with $\C^{p}$. Let $\pi: \C^{\ell} \to \C^{p}$ 
be the natural projection.
The arrangement $\A_{X}=\{H \in 
\A \ | \ H\supseteq X\}$, considered as an arrangement in 
$\C^{\ell}/X\cong \C^{p}$, determines a stratification of 
$\C^{p}$ as above.  Elements of $L(\A_{X})$ have the 
form $\pi Y=X+Y/X$ for $Y\in L(\A)$. 
Referring to Lemma \ref{sum}, one sees that the preimage of a stratum 
is a union of strata, that is, that $\pi$ is a map of stratified 
spaces, precisely when $X$ is modular. Since $\pi$ is a
linear surjection, it restricts to a submersion on each stratum. In 
order to apply the Thom Isotopy Lemma, it is necessary to extend 
$\pi$ to a proper map of stratified spaces. This step was carried 
out by L.~Paris \cite{P3}.

\begin{theorem}
There exists a stratified space $P_X$ 
containing $\C^\ell$ as an open dense subset, and an extension of 
$\pi$ to a proper stratified map $\hat{\pi}: P_X \to \C^p$.\qed
\label{paris}
\end{theorem}

The space $P_X$ is obtained by compactifying the 
fibers of $\pi$,	i.e., the parallel translates of $X$,
via projective completion, so that $P_X$ is diffeomorphic 
to $\P(\C^q)\times \C^p$, where $q=\ell-p=\dim(X)$. 
This can be viewed as a parametrized version of R.~Randell's 
construction in his proof of the lattice	isotopy	theorem \cite{R1}. 
The stratification of $\C^\ell$ is extended to a
stratification of $P_X$ by adjoining closed strata formed by	intersecting the 
closures	of the $S_Y$ in	$P_X$ with $(\P(\C^q) -	\C^q) \times \C^p$.	
These new strata	have the form 
$S_Y^\infty \times \C^p$, 
for $Y\geq X$, where	$S_Y^\infty	= (S_Y\cap X) \cap (\P(\C^q)-\C^q)$. 
The map $\hat{\pi}$ is projection on the second factor. 

Let $M(\A)$ and $M(\A_{X})$ denote the complements of $\bigcup \A$ 
and $\bigcup \A_{X}$ in $\C^{\ell}$ and $\C^{p}$ respectively.

\begin{corollary} The map $\pi|_{M(\A)}: M(\A)\to M(\A_{X})$ is a 
fiber bundle projection.
\label{bundle}
\end{corollary}

\begin{proof} The complement $M(\A)$ coincides with the open 
stratum $S_{0_{L}}$ of $\C^{\ell} \subset P_X$. So Theorem \ref{paris} implies 
that the restriction of $\widehat{\pi}$ to 
$M(\A)$ is a fiber bundle projection, by the Thom Isotopy Lemma 
\cite{GorMac,Math,Thom}.
\end{proof}

\medskip
We proceed to generalize the properties of strictly linear 
fibrations \cite{FR2}, where $X$ is a modular copoint, to general modular 
fibrations.
Henceforth let $X$ be a modular flat of $L(\A)$, and let us denote the bundle 
projection $\pi|_{M(\A)}$ by $\pi_{X}$.

We say \A\ is a {\em $K(\pi,1)$ arrangement} if $M(\A)$ is an 
aspherical space.

\begin{corollary} If $\A_{X}$ and the coned fiber arrangement $c\A_{\bar{v}}$ are 
$K(\pi,1)$ arrangements, then \A\ is a $K(\pi,1)$ arrangement.
\label{aspher}
\end{corollary}

\begin{proof} This follows immediately from the long exact homotopy sequence 
of the fibration $\pi_X$.
\end{proof}

\begin{remark}
In  case $X$ is a modular copoint, the monodromy of $\pi_{X}$ induces 
a homomorphism from $\pi_1(M(\A_X))$ to $P_{n}$, the pure braid group 
on $n=|\A-\A_{X}|$ strands, 
which we call the {\em braid monodromy homomorphism} after its 
similarity to the Moishezon construction. See \cite{CS2}.
For a 
modular flat $X$ of arbitrary rank, the pure braid group is replaced by 
the fundamental group of a certain subvariety of the Grassmanian, a 
{\em matroid stratum} defined as follows. If $P\in \G_{\ell}(\C^{n})$ is a 
point of the Grassmannian of $\ell$-planes in $\C^{n}$, then $P$ 
determines a vector configuration in $\C^{\ell}$, unique up to linear 
change of coordinates, obtained by 
projecting the standard basis vectors of $\C^{n}$ onto $P$ 
\cite{GGMS}. 
Let $G_{P}$ denote the linear matroid realized by this configuration; 
$G_P$ is independent of the choice of basis in $P$.
The {\em matroid stratum} of an arbitrary matroid $G$ is the subset 
$\Gamma(G)$ of $\G_{\ell}(\C^{n})$ given by 
$$\Gamma(G)=\{P\in \G_{\ell}(\C^{n}) \ |  \ G_{P}=G\}.$$

An ordered arrangement $\A=\{H_1,\ldots,H_n\}$ of rank 
$\ell$, with 
specified defining forms $\{\phi_1,\ldots,\phi_n\}$, determines a 
point $P\in \G_{\ell}(\C^{n})$ given by the image of 
$(\phi_1,\ldots, \phi_n): \C^\ell \to \C^n.$ The original 
arrangement \A\ is isomorphic to the arrangement in $P$ formed by the 
intersections of $P$ with the coordinate hyperplanes in $\C^n$, and 
the point $P$ lies in $\Gamma(G(\A))$. See \cite{F2}.

The monodromy of the stratified map $\pi$ induces a 
homomorphism $$\pi_{1}(M(\A_{X})) \to \pi_{1}(\Gamma(\bar{T}_{X}(G))).$$
Indeed, a path $\{\bar{v}_t\}_{t\in [0,1]}$ in the base space $M(\A_X)$ 
determines a one-parameter 
family of (coned) fiber arrangements $c\A_{\bar{v}_t}$, equipped with 
ordered sets of defining forms inherited from a fixed set of defining 
forms for \A. By Theorem \ref{lattice} and the construction above, 
this defines a path in the matroid stratum $\Gamma(\bar{T}_{X}(G))$. 
From this one easily obtains the monodromy homomorphism described above.

This construction does indeed generalize the corank one case. 
For in this case $\bar{T}_{X}(G)$ is 
a uniform matroid 
of rank two, $\Gamma(\bar{T}_{X}(G))$ is configuration space, and 
$\pi_1(\Gamma(\bar{T}_{X}(G)))$ is the pure braid group.\qed
\label{monodromy}
\end{remark}

\begin{theorem} The monodromy action of $\pi_{1}(M(\A_{X}))$ on the 
fiber $M(\A_{\bar{v}})$ is cohomologically trivial.
\label{trivial}
\end{theorem}

\begin{proof}  Since the fiber $M(\A_{\bar{v}})$ is the complement of an
arrangement, the cohomology of $M(\A_{\bar{v}})$ is free abelian, and 
is generated by $H^{1}(M(\A_{\bar{v}}))$. 
First of all we argue that the monodromy action on 
$H^{1}(M(\A_{\bar{v}}))$ is trivial, by the same reasoning as in the 
corank-one case \cite{FR2}. The group $H^{1}(M(\A_{\bar{v}}))$ 
has a free basis consisting of elements dual to the hyperplanes of 
$\A_{\bar{v}}$. Using this basis, it is clear that elements of 
$H^{1}(M(\A_{\bar{v}}))$ are uniquely determined by their 
linking numbers with the hyperplanes of $\A_{\bar{v}}$. By naturality, 
these linking numbers agree with linking numbers in $\C^\ell$ with the
hyperplanes of $\A - \A_X$. 
Since these linking numbers take values in a discrete space, and vary 
continuously, they remain locally constant under translation of the 
fiber, and thus are globally constant under translation around a 
loop in the base. This proves triviality in degree one. Since 
$H^{*}(M(\A_{\bar{v}}))$ is generated by $H^{1}(M(\A_{\bar{v}}))$, 
and the monodromy action respects cup products, it 
follows that the monodromy acts trivially on $H^{*}(M(\A_{\bar{v}}))$.
\end{proof}

A {\em rational $K(\pi,1)$ arrangement} is an arrangement whose 
complement has aspherical rational completion. 
See \cite{F4,FR1,OT,FR4} for the precise definition and basic 
properties. We point out that this property seems to bear little 
relationship to the notion of $K(\pi,1)$ arrangement; the terminology 
arises naturally in the context of simply-connected spaces. 

\begin{corollary} If $\A_{X}$ and $c\A_{\bar{v}}$ are rational 
$K(\pi,1)$ arrangements, then \A\ is a rational $K(\pi,1)$ arrangement.
\end{corollary}

\begin{proof} The argument is the same as in the corank-one case 
\cite{F4}. 
Because the monodromy action is trivial, hence 
nilpotent, on the cohomology of the fiber, the map $\pi_X$ induces a 
fibration of the rational completion of $M(\A)$ over that of 
$M(\A_X)$, with fiber the rational completion of $M(\A_{\bar{v}})$. 
Since $M(c\A_{\bar{v}})\cong \C^* \times M(\A_{\bar{v}})$, the 
hypothesis implies that the rational completion of $M(\A_{\bar{v}})$ 
is aspherical.
The assertion then follows from the homotopy sequence of this fibration.
\end{proof}

At this point the only known examples of rational $K(\pi,1)$ 
arrangements are supersolvable.
If $\A_X$ and $c\A_{\bar{v}}$ are supersolvable, 
then \A\ is also supersolvable \cite{Ox}.
So the preceding corollary does not provide new examples of rational 
$K(\pi,1)$ arrangements.

\medskip
The {\em Poincar\'e series} of a topological space $M$ is 
$$P(M,t)=\sum_{n\geq 0} \dim_\Q H^n(M,\Q).$$
For a complex 
arrangement \A, a famous result of Orlik and Solomon \cite{OS1} 
relates the Poincar\'e series $P(M(\A),t)$ to the characteristic 
polynomial $\chi(G(\A),t)$ of the underlying matroid $G(\A)$. 
Specifically, $$P(M(\A),t)=t^r\chi(G(\A),-t^{-1}),$$ where 
$r$ is the rank of $G(\A)$.

For a modular flat $X$, R.~Stanley proved in \cite{Stan3} that the characteristic polynomial 
of the $G(\A_X)$ divides that of $G(\A)$ over the integers. In 
\cite{Bry2}, T.~Brylawski identified the quotient as the 
characteristic polynomial of the complete principal truncation 
$\bar{T}_X(G)$, 
divided by $(t-1)$. The decone 
operation on arrangements has the effect on Poincar\'e polynomials of 
dividing by $(1+t)$. Using Theorem \ref{lattice} and the identity relating the characteristic 
polynomial of $G(\A)$ to the Poincar\'e polynomial of $M(\A)$, we may 
restate the Stanley and Brylawski results as follows.

\begin{theorem} If $X$ is a modular flat of $G$, then 
$$P(M(\A),t)=P(M(\A_X),t)P(M(\A_{\bar{v}}),t).$$\qed
\label{stanley}
\end{theorem}

\begin{corollary} The Leray-Serre spectral sequence of $\pi_X$ 
satisfies $$E_2^{p,q}\cong H^{p}(M(\A_{X}))\otimes 
H^{q}(M(\A_{\bar{v}})),$$ and degenerates at the $E_2$ term.\qed
\label{degenerate}
\end{corollary}

\begin{proof} The first assertion follows from the triviality of the 
monodromy action established in Theorem \ref{trivial}. 
The second is a consequence of the 
factorization identity among the Poincar\'e series. Indeed, according to 
\cite[Theorem 11.3]{Hu}, the formula of Theorem \ref{stanley} holds 
for a general spectral sequence $E$, with a correction term 
that vanishes precisely when the differential of $E_2$ is trivial.
\end{proof}

\begin{corollary} The cohomology $H^{*}(M(\A))$ is isomorphic 
as a \Z-module to 
the tensor product $H^{*}(M(\A_{X}))\otimes 
H^{*}(M(\A_{\bar{v}}))$.\qed
\label{factor}
\end{corollary}

\begin{remark} In \cite{T2} Terao established the
tensor product factorization of Corollary \ref{factor} in terms of 
Orlik-Solomon algebras, using a direct combinatorial argument. This 
approach 
yields an alternate proof of the Stanley factorization theorem.
The degeneracy of the spectral sequence in case $X$ is 
a modular copoint is given a direct proof in \cite{FR2}, providing a 
topological proof of Terao's result in this case. 
The proof in \cite{FR2} uses the fact that the fiber $M(\A_{\bar{v}})$ has nonvanishing 
cohomology in only two different degrees, so that the spectral 
sequence results in a ``Gysin-like'' long exact sequence. The other 
ingredient is the construction of a section of the
bundle map $\pi_X$. In case $X$ is a modular flat of 
arbitrary rank, a section of $\pi_X$ is constructed by L.~Paris in 
\cite{P3}. But we see no analogue of the Gysin long exact sequence in 
the general case, and do not have a topological proof, independent of the Stanley and Brylawski 
results, of the second 
part of Corollary \ref{degenerate}.
Nevertheless, the bundle map $\pi_X$ is seen to be a  
topological realization of the combinatorial and algebraic factorizations 
arising from a modular flat.\qed
\label{gysin}
\end{remark}

Motivated by the fact that supersolvable arrangements are inductively 
free \cite{JT}, we include with this compendium of generalizations the 
following conjecture.

\begin{conjecture} If $X$ is a modular flat and both $\A_X$ and 
$c\A_{\bar{v}}$ are free arrangements, then \A\ is a free arrangement.
\end{conjecture}

\end{section}

\begin{section}{Parallel connections}
\label{parallel}

Let $G_1$ and $G_2$ be matroids on ground sets $E_1$ and $E_2$. 
Suppose $E_1\cap E_2=F$ is a flat of both $G_1$ and $G_2$, and is 
modular in $G_1$. The {\em generalized parallel connection} of $G_1$ 
and $G_2$ along $F$ is the matroid $P_F(G_1,G_2)$ on the ground set 
$E_1\cup E_2$ whose flats are those sets $Y\subseteq E_1\cup E_2$ 
for which $Y\cap E_i$ is a flat of $G_i$ for $i=1,2$. The modularity 
condition is necessary for this definition to make sense. That is, 
this collection of flats will form a geometric lattice for general $G_2$ 
if and only 
if $X$ is modular in $G_1$. Modularity of $F$ in $G_1$ implies that 
$G_2$ is modular in $P_F(G_1,G_2)$. See \cite[Section 7.6]{Wh1} and \cite{Ox} 
for details about this construction.

The rank of a flat $Y$ of $P_F(G_1,G_2)$ 
is given by $$r(Y)=r_1(Y\cap E_1)+r_2(Y\cap E_2)-r_1(Y\cap F),$$
where $r_i$ is the rank function of $G_i, \ i=1,2$.
In particular, the rank of $P_F(G_1,G_2)$ is equal to 
$r(G_1)+r(G_2)-r(F).$
The rank formula indicates that $P_F(G_1,G_2)$ is the ``free'' 
sum of $G_1$ and $G_2$ amalgamated along their common flat $F$. 
Indeed, $P_F(G_1,G_2)$ is a pushout of the inclusion maps $F 
\hookrightarrow G_i, \ i=1,2$ in the category of matroids and injective strong 
maps. 

In case $F$ is a point, automatically modular in $G_1$, the matroid 
$P_F(G_1,G_2)$ is called a {\em parallel connection} of $G_1$ and 
$G_2$, studied in connection with complex hyperplane arrangements in 
\cite{EF}.

Now suppose $\A_1$ and $\A_2$ are hyperplane 
arrangements realizing $G_1$ and $G_2$ in $\C^r$ and $\C^s$ respectively. 
Then there is an arrangement \A\ realizing $P_F(G_1,G_2)$, provided 
there is a linear isomorphism between the subarrangements of $\A_1$ and $\A_2$ 
corresponding to the common flat $F$. To carry 
out the construction, let us be more precise about the realizations 
$\A_1$ and $\A_2$.

Suppose the flat $F$ has rank $p$ (in both $G_1$ and $G_2$). 
Let $X_1$ denote the corresponding element of 
intersection lattice $L(\A_1)$. Thus $X_1$ is a linear subspace of 
$\C^r$, and we may identify $(\A_1)_{X_1}$ with $F$. We may assume $X_1= \C^{r-p} \times 
\{0\} \subseteq \C^r$. Then the defining equations of the 
hyperplanes in $(\A_1)_{X_1}\subseteq \A_1$ involve only the last $p$ 
coordinates in $\C^r$. Assume that the same defining forms, expressed 
in terms of the first $p$ 
coordinates of $\C^s$, give the defining equations for 
hyperplanes  of $(\A_2)_{X_2}\subseteq \A_2$, where $X_2\in L(\A_2)$ 
corresponds to the flat $F$ of $G_2$. Then we may define an arrangement \A\ in 
$\C^\ell$, with $\ell=r+s-p$, as follows. Identify $\C^\ell$ with 
$\C^{r-p}\times \C^p \times \C^{s-p}$. By pulling back the defining 
equations via projection of coordinates, the arrangements $\A_1$ and 
$\A_2$ naturally embed in $\C^r \times \C^{s-p} = \C^\ell$ 
and $\C^{r-p}\times \C^s = \C^\ell$ respectively. Then let 
\A\ be the union of $\A_1 - (\A_1)_{X_1}$ and $\A_2$ in $\C^\ell$.

\begin{theorem} {\em\rmfamily\cite[Prop. 7.6.11]{Wh1}} The arrangement \A\ is a 
realization of the generalized parallel connection $P_F(G_1,G_2)$.
\qed
\end{theorem}

Let $X\in L(\A)$ correspond to the flat $F$ of $P_F(G_1,G_2)$.
By modularity of $F$ in $G_1$ and of $G_2$ in $P_F(G_1,G_2)$, 
the results of Section \ref{fibration} yield bundle maps $M(\A_1)\to 
M((\A_1)_{X_1})$ and $M(\A)\to M(\A_2)$. We consider $(\A_1)_{X_1}$ to be an 
arrangement in $\C^r/X_1 \cong \C^s/X_2 \cong \C^p$. Then there is a 
projection $M(\A_2)\to M((\A_1)_{X_1})$. This projection is just the 
inclusion of $M(\A_2)$ into the complement of the subarrangement 
$\{H \in \A_2 \ | \ H\supseteq X_2\}$ of $\A_2$, followed by a homotopy 
equivalence. 

\begin{theorem} The fiber bundle $M(\A)\to M(\A_2)$ is the pullback 
of the bundle $M(\A_1)\to M((\A_1)_{X_1})$ along the projection
$M(\A_2)\to M((\A_1)_{X_1})$. That is, 

\begin{equation*}
\begin{CD}
M(\A)@>>>M(\A_1)\\
@VVV @VVV\\
M(\A_2)@>>>M((\A_1)_{X_1})
\end{CD}
\end{equation*}

is a pullback diagram. 
\label{pullback}
\end{theorem}

\begin{proof} For these special realizations, the bundle map $M(\A_1)\to 
M((\A_1)_{X_1})$ 
is the restriction of the projection $\pi_1: \C^r\to \C^p$ onto 
the last $p$ coordinates. Similarly, the map $M(\A)\to M(\A_2)$ is 
the restriction of the projection $\pi: \C^\ell \to \C^s$ onto the last 
$s$ coordinates. The map $M(\A_2)\to M((\A_1)_{X_1})$ 
can be identified with the restriction of the projection $\pi_2: \C^s \to 
C^s/X_2 \cong \C^p$ onto the first $p$ coordinates.

By definition, the total space of the pullback of $M(\A_1)\to M((\A_1)_{X_1})$ 
along the projection $M(\A_2)\to M(\A_X)$ is the set 
of pairs $(x,v)\in \C^r \times \C^s$ such that $x\in M(\A_1), \ v\in 
M(\A_2)$, and $\pi_1(x)=\pi_2(v)$ in $M((\A_1)_{X_1})\subseteq C^p$. But this 
means that the last $p$ components of $x$ match the first $p$ 
components of $v$. Then each such $(x,v)$ corresponds to a unique point of 
$\C^{r+s-p} =\C^\ell$ which, by the first two conditions, lies in 
$M(\A)$. Under this identification, the projection $(x,v)\mapsto v$ 
coincides with $\pi$. This identifies $\pi|_{M(\A)}$ with the pullback 
of $\pi_1|_{M(\A_1)}$, as claimed.
\end{proof}

\begin{corollary} If $\A$ is a realization of the parallel connection 
of $\A_1$ and $\A_2$, then $M(\A)$ is a trivial bundle over $M(\A_2)$ 
with fiber $M(d\A_1)$. In particular, $M(\A)\cong M(d\A_1)\times 
M(\A_2)$.
\end{corollary}

\begin{proof} In case $X$ is a point, then $X$ is modular in $\A_1$, and the 
modular fibration $M(\A_1)\to M(\A_X) = \C^*$ is a trivial bundle with 
fiber $d\A_1$, by \cite[Proposition 5.1]{OT}. The pullback of a trivial 
bundle is trivial.
\end{proof} 

This corollary clarifies the main construction of \cite{EF}, which 
essentially established the diffeomorphism noted above. This argument 
shows in an alternate way that the diffeomorphisms among arrangements 
with different underlying matroids, constructed in \cite{EF}, are all 
consequences of the triviality of the restriction of the Hopf bundle.
\end{section}

\begin{section}{Examples}
\label{examples}

Corollary \ref{aspher} of Section \ref{fibration} can be used to 
identify $K(\pi,1)$ arrangements of high rank, at least when the base 
arrangement $\A_X$ and (coned) fiber arrangement $c\A_{\bar{v}}$ are 
tractable. This will be the case, for instance, when $X$ is a 
modular coline, for then $c\A_{\bar{v}}$ will have rank three. In this 
section we present new families  of examples of $K(\pi,1)$ 
arrangements. Our results give some support for the conjecture 
\cite{FR4}, which was based 
primarily on rank-three phenomena, that factored arrangements of 
arbitrary rank are $K(\pi,1)$. We also exhibit an 
interesting example with two different modular colines, allowing 
us to conclude the nontrivial result that one of the fiber 
arrangements is not $K(\pi,1)$.

Let $\B_\ell$ denote the arrangement of reflecting hyperplanes in the 
Weyl group of type $B_\ell$. Thus $\B_\ell$ consists of the hyperplanes 
$$H_{ij}=\{x\in \C^\ell \ | \ x_i=x_j\}, \ \text{for} \ 1\leq i<j\leq 
\ell, $$
$$\bar{H}_{ij}=\{x\in \C^\ell \ | \ x_i=-x_j\} \ \text{for} \ 1\leq i<j\leq 
\ell, \ \text{and}$$
$$H_i=\{x\in \C^\ell \ | \ x_i=0\}, \ \text{for} \ 1\leq i\leq \ell.$$
Let $\A_{\ell-1}$ denote the braid arrangement, consisting of the 
hyperplanes $H_{ij}$ above, for $1\leq i<j\leq \ell$.

In \cite{ER}, P.~Edelman and V.~Reiner used graphs to parametrize 
subarrangements of $\B_\ell$ containing $\A_{\ell-1}$,
developing a calculus for combinatorial invariants of the 
arrangements in terms of the graphs. We find among these arrangements 
those which are not supersolvable, but have modular colines, for which the 
fiber arrangements 
are demonstrably $K(\pi,1)$. These examples coincide in large part 
with the arrangements 
between $\A_{\ell-1}$ and $\B_\ell$ which are factored, classified 
by D.~Bailey in \cite{Bail}.

Let $\Gamma$ be a graph with vertex set 
$\{1,\ldots,\ell\}$, possibly with loops, but without multiple edges. 
Let $\edge(\Gamma)$ and $\lop(\Gamma)$ denote the sets of edges and 
loops of $\Gamma$, respectively.
Let $\A_\Gamma$ be the arrangement defined by 
$$\A_\Gamma=\A_{\ell-1} \ \cup \ \{\bar{H}_{ij} \ | \ ij\in 
\edge(\Gamma)\} \ \cup \ \{H_i \ | \ i\in \lop(\Gamma)\}.$$

The following results are proved in \cite{ER}. The notion
of free arrangement plays little role in what follows; see \cite{OT} for a 
precise definition. A graph is {\em threshold} if it is built up by 
successively adjoining isolated and/or cone vertices, the latter 
being vertices which are 
adjacent to all preceding vertices. 

\begin{theorem} The arrangement $\A_\Gamma$ is free if and only if
\begin{enumerate}
\item $\Gamma$ is a threshold graph, and
\item $i\in \lop(\Gamma)$ and $\deg(j)>\deg(i)$ implies 
$j\in\lop(\Gamma)$. \qed
\end{enumerate}
\label{free}
\end{theorem}

An edge $ij\in edge(\Gamma)$ is {\em loopless} if neither $i$ nor 
$j$ lies in $\lop(\Gamma)$.

\begin{theorem} Suppose $\A_\Gamma$ is free and $\Gamma$ has no loopless edges. Then 
$\A_\Gamma$ is supersolvable.\qed
\label{ss}
\end{theorem}

There are two families of exceptional graphs $\Gamma$ with loopless 
edges such that $\A_\Gamma$ 
is supersolvable; see \cite{ER}.
Of course, any such arrangement $\A_\Gamma$ is $K(\pi,1)$.

\medskip
Roughly speaking, an arrangement is {\em factored} \cite{J2,FJ,T4,OT} if 
the cohomology of the complement
is isomorphic as a \Z -module to the the tensor product of algebras with 
trivial multiplication
generated by sets of hyperplanes of \A. 
For instance, Corollary \ref{factor} implies that supersolvable 
arrangements are factored.
Such factorizations correspond to partitions of \A, properties of 
which are analyzed in \cite{FJ,T4,JP}. D.~Bailey \cite{Bail} identified 
those arrangements $\A_\Gamma$ which are factored.

\begin{theorem} {\em\rmfamily \cite{Bail}} Suppose $\A_\Gamma$ is free and
$\Gamma$ has at most one loopless edge. Then 
$\A_\Gamma$ is a factored arrangement. \qed
\label{bailey}
\end{theorem}
Again there are two families of exceptional graphs with more than one 
loopless edge for which $\A_\Gamma$ is factored \cite{Bail}.

\medskip
We establish criteria for $\A_\Gamma$ to have modular copoints or 
colines determined by coordinate subspaces $X$. The assertions below are easy 
to prove using Remark \ref{converse}, by showing that the lattices of 
the coned fiber 
arrangements remain constant over $M(\A_X)$. (For an example, see 
the proof of Theorem \ref{new}.) 
Let $\Gamma'$ and $\Gamma''$ be the vertex-induced subgraphs of $\Gamma$ 
on vertices $\{2,\ldots,\ell\}$ and $\{3,\ldots,\ell\}$ respectively. 
Then $\A_{\Gamma'}$ and $\A_{\Gamma''}$ are flats of 
$G(\A_\Gamma)$ of corank one and two corresponding to the linear subspaces 
$x_2=\cdots=x_\ell=0$ and $x_3=\cdots =x_\ell=0$ respectively.

\begin{lemma} The flat $\A_{\Gamma'}$ is a modular copoint of 
$G(\A_\Gamma)$ if and only if 
\begin{enumerate}
\item $1\in \lop(\Gamma)$ implies $\lop(\Gamma)=\{1,\ldots,\ell\}$, 
and
\item $1j\in \edge(\Gamma)$ implies $j\in \lop(\Gamma)$ and 
$j$ is adjacent \\ to every vertex
of $\Gamma$. \qed 
\end{enumerate}
\label{copoint}
\end{lemma}
In particular, an isolated vertex of $\Gamma$ corresponds to a 
modular copoint.

\begin{lemma} The flat $\A_{\Gamma''}$ is a modular coline of 
$G(\A_\Gamma)$ if and only if
\begin{enumerate}
\item  $i\in \lop(\Gamma)$, for $i=1$ or $2$ implies 
$\lop(\Gamma)\supseteq\{3,\ldots,\ell\}$,
\item $ij\in \edge(\Gamma)$ for $i=1$ or $2$ implies
$j\in\lop(\Gamma)$ and $j$ is adjacent to every vertex $k$ for $k\geq 
3$, and
\item $12\in\edge(\Gamma)$ implies $\lop(\Gamma)\supseteq 
\{3,\ldots,\ell\}$. \qed \\
\end{enumerate}
\label{coline}
\end{lemma}

The modular fibration corresponding to a modular copoint has as coned 
fiber arrangement $c\A_{\bar{v}}$ a central arrangement of rank two, 
which is $K(\pi,1)$, so $\A_\Gamma$ is $K(\pi,1)$ if and only if 
$\A_{\Gamma'}$ is $K(\pi,1)$. This is a ``strictly linear fibration" 
as studied in \cite{FR2,T2}.
 
Factored arrangements of rank three were shown to be $K(\pi,1)$ in 
\cite{P2}. Supersolvable arrangements of arbitrary rank are factored, 
and are $K(\pi,1)$. In \cite{FR4} we conjecture that factored 
arrangements of arbitrary rank are $K(\pi,1)$. The next result 
provides more support for this conjecture. 
The arrangements of this theorem are not supersolvable, by Theorem \ref{ss}, but are 
factored, by Theorem \ref{bailey}. In fact, by an argument in 
\cite{Bail}, the examples described below are the only factored 
non-supersolvable arrangements $\A_\Gamma$ which have only two non-loop vertices.

\begin{theorem}	Suppose	$\Gamma$ is the complete graph on $\ell$ 
vertices, and $\lop(\Gamma)=\{3,\ldots,\ell\}$. 
Then $\A_\Gamma$ is a $K(\pi,1)$ arrangement.
\label{new} 
\end{theorem}

\begin{proof}  Let $\A=\A_\Gamma$. Let $\Gamma''$ denote the vertex-induced subgraph of 
$\Gamma$ on vertices $\{3,\ldots,\ell\}$. Then $\A_{\Gamma''}$ is a modular coline of 
$G(\A)$ by Lemma \ref{coline}, and $\A_{\Gamma''}$ 
is supersolvable, hence $K(\pi,1)$, by Theorem \ref{ss}. 
Let $X=\bigcap \A_{\Gamma''}$, and let $\bar{v}=(v_3,\ldots, v_\ell)\in M(\A_X).$ 
Then $v_i\not = \pm v_j$ for $3\leq i<j\leq \ell$, and $v_j\not=0$ for $3\leq j\leq \ell$. 
The fiber arrangement $\A_{\bar{v}}$ is the affine arrangement in $\C^2$ consisting of the lines 
$$x_1=\pm x_2, \ x_1=\pm v_j, \ \text{and} \ x_2=\pm v_j$$ for $3\leq j\leq \ell.$ 

\begin{figure}[h]
\begin{center}
\epsfig{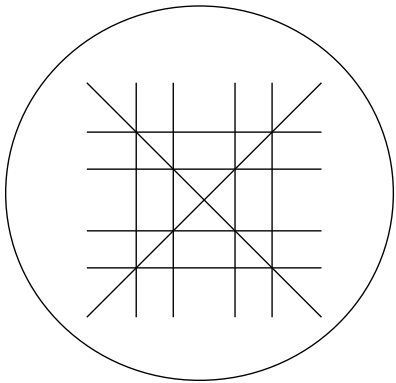}
\caption{The coned fiber arrangement, for $\ell=4$}
\label{fiber}
\end{center}
\end{figure}

The coned fiber arrangement 
$c\A_{\bar{v}}$ pictured in Figure~\vref{fiber} is a $K(\pi,1)$ arrangement. 
Indeed $c\A_{\bar{v}}$ is precisely 
Example 3.13 of \cite{F1}. Alternatively, $c\A_{\bar{v}}$ is a 
factored arrangement of rank three, hence $K(\pi,1)$ by \cite{P2}. We 
conclude by Corollary \ref{aspher} that \A\ is $K(\pi,1)$.
\end{proof}

\begin{remark} The conclusion of the theorem also holds if 
$\lop(\Gamma)=\emptyset$, for then $\A_\Gamma$ is a Coxeter arrangement 
of type $D_\ell$, which is simplicial.\qed
\end{remark}

One can use Theorem \ref{new} and Lemma \ref{copoint} 
to build other examples of non-supersolvable (and non-simplicial) 
$K(\pi,1)$ arrangements of high rank, by successively adding vertices 
satisfying \ref{copoint} to the graphs of Theorem \ref{new}. 
See Figure \vref{graph1}. 
The existence of loops in $\Gamma$ is essential: the same 
construction with the loopless complete graph (the $D_\ell$ 
arrangement) allows only the addition of isolated vertices.

There is one rank-three arrangement $\A_\Gamma$, a realization of the 
non-Fano matroid, which is not 
supersolvable, but is simplicial, hence $K(\pi,1)$. 
This arrangement can also be used with Theorem \ref{copoint}
to construct non-supersolvable $K(\pi,1)$ arrangements. 
This construction is illustrated in Figure \vref{graph3}.

\begin{example} 
We exhibit in Figures \ref{graph1} and \vref{graph3}
some other graphs $\Gamma$ for which $\A_\Gamma$ 
is $K(\pi,1)$, using the constructions of the preceding paragraphs.
These arrangements are factored, as is every arrangement arising from these 
constructions, by Theorem \ref{bailey}.  
\qed
\end{example}

\begin{figure}[h]
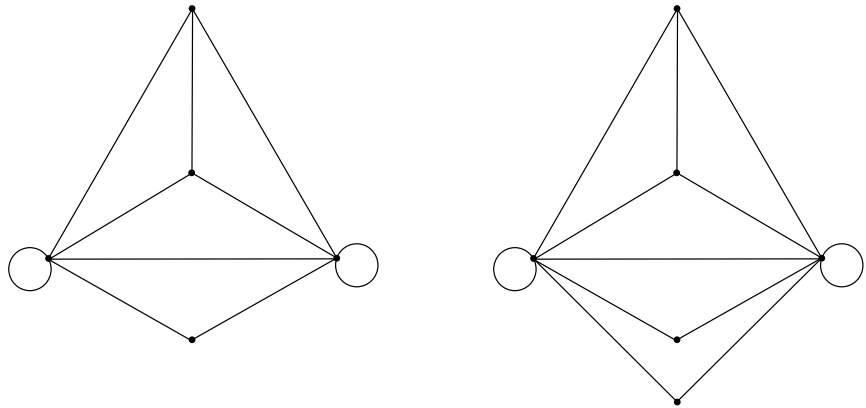

\begin{center}
\mbox{\subfigure{\epsfig{file=fig2a.epsf, height=4.5cm}}
\qquad \qquad \subfigure{\epsfig{file=fig2b.epsf, height=4.5cm}}}
\vspace*{.5cm}
\caption{Extensions of Theorem \ref{new}}\label{graph1}
\end{center}
\end{figure}

\begin{figure}[h]
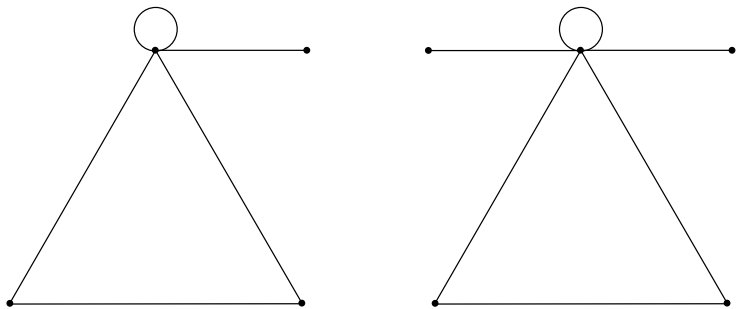

\begin{center}
\mbox{\subfigure{\epsfig{file=fig3a.epsf, height=4cm}} \qquad \qquad
\subfigure{\epsfig{file=fig3b.epsf, height=4cm}}}
\caption{Graphs of $K(\pi,1)$ extensions of a non-Fano arrangement}\label{graph3}
\end{center}
\end{figure}

We close with another interesting example from the class of 
``$A$-$B$ arrangements.'' 

\begin{example} Let $\Gamma$ be the graph with vertex set 
$\{1,2,3,4\}$, edges $12, \ 13,$ and $24$, and 
loops at vertices $1,2, \ \text{and} \ 3$. Let $\Gamma_0''$ and 
$\Gamma_1''$ be the vertex-induced subgraphs of $\Gamma$ with vertex 
sets $\{1,2\}$ and $\{1,3\}$ respectively. Then both $\A_{\Gamma_0''}$ 
and $\A_{\Gamma_1''}$ are modular flats of $G(\A_\Gamma)$. The 
respective fiber arrangements are illustrated in Figure~\vref{not}. 

\begin{figure}[h]
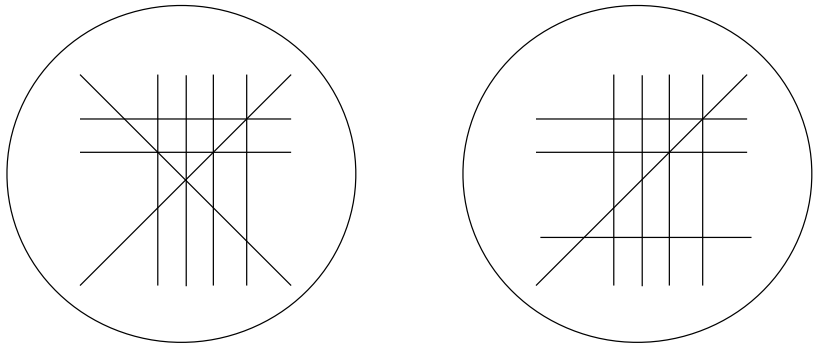

\begin{center}
\mbox{\subfigure{\epsfig{file=fig4a.epsf, height=4.5cm}}\qquad \qquad 
\subfigure{\epsfig{file=fig4b.epsf, height=4.5cm}}}
\caption{Two non-$K(\pi,1)$ arrangements}
\label{not}
\end{center}
\end{figure}

The arrangement on the right, the (coned) fiber arrangement of $M(\A_\Gamma)\to 
M(\A_{\Gamma_1''})$, is not $K(\pi,1)$ because it has a ``simple 
triangle'' \cite[Corollary 3.3]{FR1}. Now $\A_{\Gamma_1''}$ is a $K(\pi,1)$ 
arrangement, being a
central arrangement of rank two. It follows that $M(\A_\Gamma)$ is 
not aspherical. But $\A_{\Gamma_0'}$ is also $K(\pi,1)$. We conclude 
that the (coned) fiber arrangement of $M(\A_\Gamma)\to 
M(\A_{\Gamma_0'})$, shown on the left, cannot be $K(\pi,1)$. This is the only argument we 
know of to show this arrangement is not $K(\pi,1)$.
\qed
\end{example}

The research presented here, and our general interest in $K(\pi,1)$ 
arrangements, was motivated in part by a suggestion of G.~Ziegler of a 
straightforward construction of rank-three 
arrangements with the same underlying matroid but homotopy inequivalent complements.
The argument avoids fundamental group computations, but relies 
on the existence of high-rank $K(\pi,1)$ 
arrangements with certain properties, whose existence has not yet been 
shown. Here is the precise problem, to which the methods of this paper 
may apply.

\begin{problem} (Ziegler) Find $K(\pi,1)$ arrangements 
whose matroids have the same flats of ranks one and two but 
have different characteristic polynomials.
\end{problem}

\end{section}

\begin{ack} We are grateful to Joseph Kung, Luis Paris, Hiroaki Terao and 
Vic Reiner and David Bailey, for 
helpful conversations and correspondence, and to Terence Blows, who 
organized the {\em REU} program at Northern Arizona University, where 
this project was begun.
\end{ack}


\makeatletter
\bigskip
{\small
Department of Mathematics and Statistics \\
Northern Arizona University \\
Flagstaff, AZ 86011-5717 \\
michael.falk@nau.edu 
\\ 
\\
Department of Mathematics \\
Harvard University \\
Cambridge, MA  \\
Proudf@fas.harvard.edu
\makeatother

\end{document}